\input amstex
\documentstyle{amsppt}
\topmatter
%\UDclass   \endUDclass
\title
About integrability of almost complex structures on strictly Nearly K\"{a}hler 6-manifolds
\endtitle
\rightheadtext{Almost complex structures on SNK 6-manifolds}
\thanks
The work was partially supported by RFBR, 12-01-00873-a and by
Russian President Grant supporting scientific schools
SS-544.2012.1
\endthanks
\author    Natalia Daurtseva  \endauthor
%\xauthor   _.~_.~"ами<ия\endxauthor
%\date     \enddate
%\date      _кR-чат_<ь-ыc вариа-т  _.\enddate
\address Russia, Kemerovo State University
\endaddress
%\address   _.~\endaddress
\email  natali0112\@ngs.ru \endemail
\abstract We show that any almost complex structure, positively tamed with $\omega$ on nearly K\"{a}hler 6-manifold
$(M,g,J,\omega)$ is not integrable.
\endabstract
\endtopmatter
\document
{\bf 1. Introduction.} Let $(M,g,J,\omega)$ is an almost Hermitian manifold, where
$g$ is a Riemannian metric, $J$ an almost complex structure compatible with
$g$ and $\omega$ is K\"{a}hler form, $\omega(X,Y)=g(JX,Y)$.
Nearly K\"{a}hler manifold is an almost
Hermitian manifold $(M,$ $g,$ $J$, $\omega)$ with the property
that $(\nabla_XJ)X=0$ for all tangent vectors $X$, where $\nabla$
denotes the Levi-Civita connection of $g$. If $\nabla_XJ\neq 0$
for any non-zero vector field $X$, $(M,g,J,\omega)$ is called strictly
nearly K\"{a}hler. Nearly K\"{a}hler geometry comes from the
concept of weak holonomy introduced by A. Gray in 1971 [1], this
geometry corresponds to weak holonomy $U(n)$. The class of nearly
K\"{a}hler manifolds appears naturally as one of the sixteen
classes of almost Hermitian manifolds described by the
Gray-Hervella classification [2]. Nagy P.-A. has proved [3] that
every compact simply connected nearly K\"{a}hler manifold $M$ is
isometric to a Riemannian product $M_1\times\dots M_k$, such that
for each $i$, $M_i$ is a nearly K\"{a}hler manifold belonging to
the following  list: K\"{a}hler manifolds, naturally reductive
3-symmetric spaces, twistor spaces over compact
quaternion-K\"{a}hler manifolds with positive scalar curvature,
and nearly K\"{a}hler 6-manifold. This is one of the reasons of
interest to the nearly K\"{a}hler 6-manifolds. In case of
dimension 6 there exists several equivalent conditions defining a
strictly nearly K\"{a}hler structures $(g,J,\omega)$ on $M$. For
example, conditions

(i) $(\nabla_XJ)Y$ is skew-symmetric with respect to $X,Y$ and
non-zero;

(ii) The form $\nabla\omega$ is non-zero, totally skew-symmetric
and $\nabla_X\omega=\frac13\iota_Xd\omega$, $\forall
X\in\Gamma(TM)$;

(iii) The structure group of $M$ admits a reduction to $SU(3)$, that is, there is (3,0)-form
$\Omega$ with $|\Omega|=1$, and
$$
d\omega=3\lambda \,Re\,\Omega,\ d\,Im \,
\Omega=-2\lambda\omega^2
$$
where $\lambda$ is a non-zero real constant, are equivalent and
define strictly nearly K\"{a}hler manifold (see [4]).

By (iii), for strictly nearly K\"{a}hler structure $(g,J,\omega)$ the $d\omega\neq 0$
and of type (3,0)+(0,3), so $J$ is not integrable. It is
natural to ask about the possibility of "neighborhood" of nearly
K\"{a}hler structure $J$ with integrable one.
%As fundamental form $d\omega$ is of type (3,0)+(0,3), with respect to nearly K\"{a}hler structure $J$,
%then, it is interesting to understand, do there exists some structure $I$, compatible with $\omega$, but integrable.
%In this note give answer on this question.
It is known [5] that the almost complex structure, underlying a strictly nearly
K\"{a}hler 6-manifold cannot be compatible with any symplectic form, even locally.
In this article we give answer on question: do there exist integrable almost complex structure, compatible with $\omega$?

{\bf 2. The main result.} Denote by:

$A^+$ --  the space of all positive oriented almost complex structures (a.c.s.) on
$M$.

$A^+_{\omega}$ --  the set of all almost complex structures on $M$, tamed by form $\omega$.

$$
A^+_{\omega}=\{I\in A^+:\omega(IX,IY)=\omega(X,Y),  \forall X,Y\in\Gamma(TM); \omega(X,IX)>0,\forall X\neq 0\}
$$

$AO^+_g$ -- the set of all $g$-orthogonal positive oriented almost complex structures.
$$
AO^+_g=\{I\in A^+:g(IX,IY)=g(X,Y),\forall X,Y\in\Gamma(TM)\}
$$
\proclaim{Lemma}
For arbitrary a.c.s.
$I\in A^+_{\omega}$
endomorphism $J+I$ is non-degenerate.
\endproclaim

\demo{Proof}
Suppose that there exists the a.c.s.
$I\in A^+_{\omega}$,
for which
$J+I$
 is degenerate. Then
$(J+I)_x(X)=0$
for some $x\in M$, $X\in\Gamma(TM)$.
Therefore
$\omega_x(X,JX)=-\omega(X,IX)$,
what contradicts to positivity of
$\omega(X,JX)$
and
$\omega(X,IX)$.
Lemma is proved.
\enddemo

As a consequence of lemma for each a.c.s. $I\in A^+_{\omega}$ one can find endomorphism:
$$
K=(J+I)^{-1}(I-J)=(1-IJ)^{-1}(1+IJ)
$$
It is easy to show, that for $K$:

1. $KJ=-JK$;

2. $g(KX,Y)=g(X,KY)$, $\forall X,Y\in\Gamma(TM)$;

3. $1-K^2>0$

Properties 1 and 2 give additionally that
$\omega(KX,Y)=-\omega(X,KY)$. On the other hand, any endomorphism $K$, with properties 1-3 defines a.c.s.
$I=(1-K)J(1-K)^{-1}\in A^+_{\omega}$.

\proclaim{Theorem}
Any almost complex structure
$I\in A_{\omega}^+$
is non-integrable.
\endproclaim

\demo{Proof}
Let
$X\in\Gamma(TM)$
 -- is field of eigenvectors of
$1-K^2$
with eigenvalue
$\lambda$. Then
$(1-K^2)JX=J(1-K^2)X=\lambda JX$,
and
$(1-K^2)(X-iJX)=\lambda(X-iJX)$.
So one can define $g$-orthogonal basis
$\nu^1, \nu^2, \nu^3$
in the space of (1,0)-forms, which are eigenvectors of
$1-K^2$ in $T^*M\otimes\Bbb C$.

Form
$\omega(\nu^k,\nu^l)=g(J\nu^k,\nu^l)=ig(\nu^k,\nu^l)=0$,
for
$k,l=1,2,3$
so $\omega$ is of type (1,1) in the basis
$(\nu,\overline{\nu})$.

Let define forms
$\theta^k=(1-K)\nu^k$,
$k=1,2,3$
in
$T^*M\otimes\Bbb C$. One can see, that $I\theta^k=(1-K)J(1-K)^{-1}(1-K)\nu^k=i(1-K)\nu^k=i\theta^k$.
Therefore,
$\theta^1,\theta^2,\theta^3$
 -- are linear independent (1,0)-forms, with respect of $I$.

Find value of form
$\omega$
on
$\theta^k,\theta^l$
for arbitrary
$k,l=1,2,3$:
$$
\omega(\theta^k,\theta^l)=\omega((1-K)\nu^k,(1-K)\nu^l)=\omega(\nu^k,(1-K^2)\nu^l)=\lambda\omega(\nu^k,\nu^l)=0.
$$
Therefore form
$\omega$
has type (1,1) in basis
$(\theta,\overline{\theta})$.

Find $d\omega^{(3,0)}$ in $(\theta,\overline{\theta})$:
$$
d\omega=\nu^1\wedge\nu^2\wedge\nu^3+\overline{\nu}^1\wedge\overline{\nu}^2\wedge\overline{\nu}^3=
$$
$$
=(1-K)^{-1}\theta^1\wedge(1-K)^{-1}\theta^2\wedge(1-K)^{-1}\theta^3+
\overline{(1-K)}^{-1}\overline{\theta}^1\wedge\overline{(1-K)}^{-1}\overline{\theta}^2\wedge\overline{(1-K)}^{-1}\overline{\theta}^3
$$
Let's calculate the
$d\omega^{(3,0)}$.
In local frame
$(\nu,\overline{\nu})$
the matrix of operator
$1-K$
is
$\left(
\matrix
1 & \overline{V}\\
V & 1
\endmatrix\right)$,
where
$V^T=V$,
$1-V\overline{V}>0$.
It is easy to check, that
$\left(
\matrix
1 & \overline{V}\\
V & 1
\endmatrix\right)^{-1}=\left(
\matrix
\Lambda^{-1} & -\overline{V}\Lambda^{-1}\\
-V\Lambda^{-1} & \Lambda^{-1}
\endmatrix
\right)$,
where
$\Lambda$
- is diagonal matrix of eigenvalues of
$1-K^2$, corresponding to $(\nu^1,\nu^2,\nu^3)$.
Therefore
$$
d\omega^{(3,0)}=\frac{1}{\sqrt{\det(1-K^2)}}\theta^1\wedge\theta^2\wedge\theta^3\neq 0
$$
As $\omega$ is of type (1,1), then
$d\omega^{(3,0)}\neq 0$
gives
$d^{2,-1}\neq 0$
for
$I$,
and shows non-integrability of this structure.
Theorem is proved.
\enddemo

\remark{Remark} The proof of the above theorem shows that really
we needn't in nearly K\"{a}hler structure, we just use, that
$d\omega^{(3,0)+(0,3)}\neq 0$, and $d\omega^{(1,2)+(2,1)}=0$.
\endremark
{\bf 3. Almost complex structures on $S^6$.} It is known that
$S^6$ admits the set of nearly K\"{a}hler structures. All of them
are orthogonal with respect to round metric
 $g_0$,
$G_2$
invariant and form the space
$\Bbb R P^7=SO(7)/G_2$.
All
$g_0$
orthogonal almost complex structures on
$S^6$
are not integrable [6]. For Riemannian manifold
$(M,g_0)$
admitting almost complex structure [7] the space
$A^+$
is a smooth locally trivial fiber bundle over the space
$AO^+_{g_0}$,
with fiber
$A^+_{\omega_J}$
over
$J\in AO^+_{g_0}$, where $\omega_J(X,Y)=g_0(JX,Y)$.

So the above theorem let us to enlarge the number of non integrable almost complex structures by
the structures in fibers over the nearly K\"{a}hler ones.

{\it Acknowledgements.} The author is grateful to M.Verbitsky and N.Smolentsev for stimulating discussions.
%\enddocument

\enddocument